\documentclass[10pt]{amsart}
\usepackage{amssymb}
\usepackage{epsfig}
\usepackage{latexsym}
\usepackage{amsmath,amssymb,amsfonts,amsthm,graphics}
\usepackage{eucal}
\usepackage{eufrak}
\usepackage[all]{xypic}
\usepackage{xspace}

\theoremstyle{plain}
\newtheorem{theorem}{Theorem}
\newtheorem{corollary}{Corollary}

\newtheorem{lemma}{Lemma}

\theoremstyle{definition}
\newtheorem{definition}{Definition}

\theoremstyle{remark}
\newtheorem{remark}{Remark}

\numberwithin{equation}{section}

\begin{document}


\title[Crossings and nesting in tangled-diagrams]
      {Crossings and nestings in tangled-diagrams}
\author{William Y. C. Chen, Jing Qin and Christian M. Reidys$^{\,\star}$}
\address{Center for Combinatorics, LPMC-TJKLC \\
          Nankai University  \\
          Tianjin 300071\\
          P.R.~China\\
          Phone: *86-22-2350-6800\\
          Fax:   *86-22-2350-9272}
\email{reidys@nankai.edu.cn}
\thanks{}
\keywords{tangled-diagram, partition, matching, crossing, nesting, vacillating
tableau}
\date{August, 2007}
\begin{abstract}
A tangled-diagram over $[n]=\{1,\dots,n\}$ is a graph of degree less than
two whose vertices $1,\dots,n$ are arranged in a horizontal line and whose
arcs are drawn in the upper halfplane with a particular notion of crossings
and nestings.
Generalizing the construction of Chen~{\it et.al.} we prove a bijection
between generalized vacillating tableaux with less than $k$ rows and
$k$-noncrossing tangled-diagrams and study their crossings and nestings.
We show that the number of $k$-noncrossing and $k$-nonnesting
tangled-diagrams are equal and enumerate tangled-diagrams.
\end{abstract}
\maketitle {{\small
}}



\section{Introduction}\label{S:intro}


The main objective of this paper is to study tangled-diagrams by
generalizing the concept of vacillating tableaux introduced by
Chen {\it et.al.}~\cite{Chen}. Tangled-diagrams are motivated from
intra-molecular interactions of RNA nucleotides as follows:
the primary sequence of an RNA molecule is the sequence of
nucleotides {\bf A}, {\bf G}, {\bf U} and {\bf C} together with the
Watson-Crick ({\bf A-U}, {\bf G-C}) and ({\bf U-G}) base pairing
rules specifying the pairs of nucleotides can potentially form
bonds. Single stranded RNA molecules form helical structures whose
bonds satisfy the above base pairing rules and which, in many cases,
determine their function. One question of central importance is now
to predict the 3D-arrangement of the nucleotides, vital for the
molecule's functionality. For this purpose it is important to
capture the sterical constraints of the base pairings, which then
allows to systematically search the
configuration space. For a particular
class of RNA structures, the pseudoknot RNA structures
\cite{Reidys:07pseu}, the notion of diagrams \cite{Stadler:99} has been
used in order to translate the biochemistry of the nucleotide
interactions \cite{Batey,Sivakova} into crossings and nestings of
arcs. A diagram is a labeled, partial one-factor graph over $[n]$,
represented as follows: all vertices are drawn in a
horizontal line and all arcs (representing interactions) are drawn
in the upper halfplane. Since {\it a priori} restricted to degree $\le 1$
diagrams do not allow to express helix-helix, loop-helix and
multiple nucleotide interactions in general \cite{Batey}. The
tangled-diagrams studied in the following are tailored to express
these interactions. A tangled-diagram is a labeled graph over the
vertices $1,\dots,n$, drawn in a horizontal line and arcs in the
upper halfplane. In general, it has isolated points and the
following types of arcs
\begin{center}
\scalebox{0.6}[0.6]{\includegraphics*[60,650][560,820]{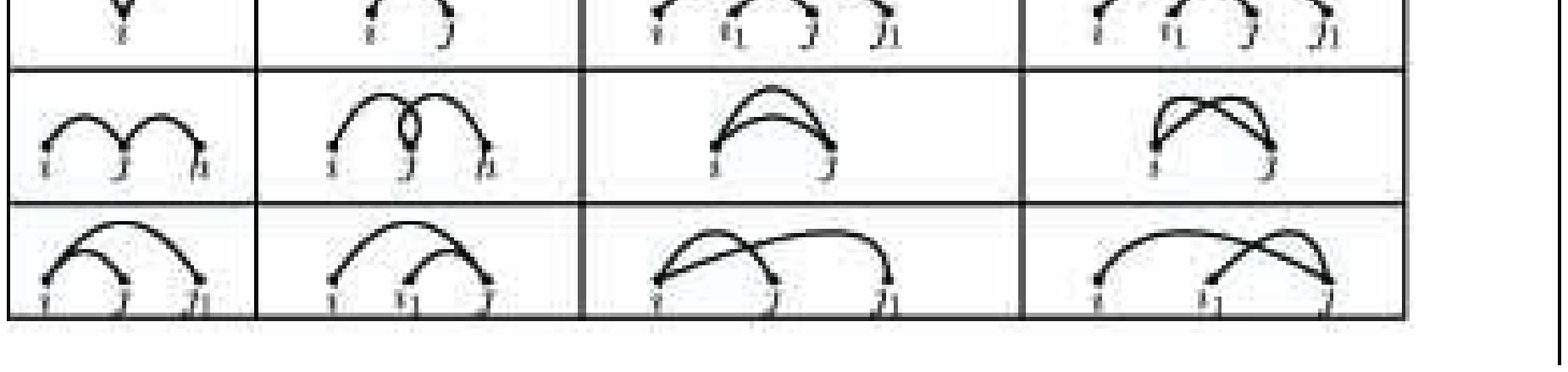}}
\end{center}
For instance \\
{\setlength{\fboxsep}{0.5pt}
\begin{center}
\scalebox{0.5}[0.5]{\includegraphics*[20,770][580,840]{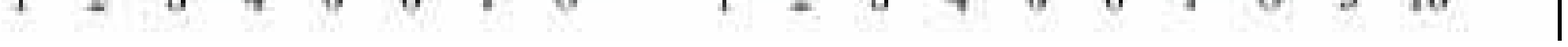}}
\end{center}}
are two tangled-diagrams and diagrams in which all vertices $j$ of
degree two are either incident to loops $(j,j)$ or crossing arcs
$(i,j)$ and $(j,h)$, where $i<j<h$ are called braids. In particular
matchings and partitions are tangled-diagrams. A matching over the
set $[2n]=\{1,2,\dots,2n\}$ is just a $1$-regular tangled-diagram
and a partition corresponds to a tangled-diagram in which any vertex
of degree two, $j$, is incident to the arcs $(i,j)$ and $(j,s)$,
where $i<j<s$. For instance
\begin{center}
\scalebox{0.6}[0.6]{\includegraphics*[10,750][580,850]{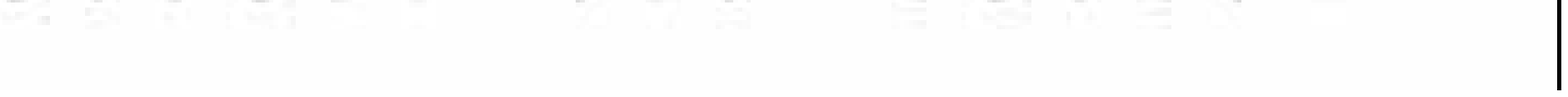}}
\end{center}
Chen {\it et al.} observed that there is a bijection between
vacillating tableaux and partitions \cite{Chen} derived from an
RSK-insertion idea due to Stanley. In addition they studied enhanced
partitions via hesitating tableaux. In the following we integrate
both frameworks by generalizing vacillating tableaux as follows: a
generalized vacillating tableaux $V_{\lambda}^{2n}$ of shape
$\lambda$ and length $2n$ is a sequence $(\lambda^{0},
\lambda^{1},\ldots,\lambda^{2n})$ of shapes such that
$\lambda^{0}=\varnothing$ and $\lambda^{2n}=\lambda,$ and
$(\lambda^{2i-1},\lambda^{2i})$ is derived from $\lambda^{2i-2}$,
for $1\le i\le n$ by an elementary move, i.e.~a step of the form
$(\varnothing,\varnothing)$: do nothing twice;
$(-\square,\varnothing)$: first remove a square then do nothing;
$(\varnothing,+\square)$: first do nothing then adding a square;
$(\pm \square,\pm \square)$: add/remove a square at the odd and even
steps, respectively. For instance the below sequence is a
generalized vacillating tableaux
\begin{center}
\scalebox{0.6}[0.6]{\includegraphics*[10,750][580,830]{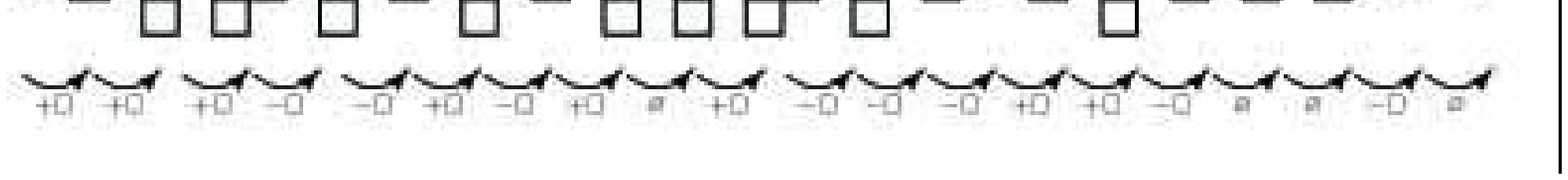}}
\end{center}
We prove a bijection between $V_\varnothing^{2n}$, referred to from
now on as simply vacillating tableaux, and tangled-diagrams over
$[n]$. We show that the notions of $k$-noncrossing tangled-diagrams
and $k$-nonnesting in tangled-diagrams are in fact dual and
enumerate $k$-noncrossing tangled-diagrams. Restricting the steps of
the vacillating tableaux we obtain three wellknown bijections, the
bijection between vacillating tableaux with elementary moves
$\{(-\square,\varnothing),(\varnothing,+\square)\}$ and matchings
\cite{Chen}, the bijection between the vacillating tableaux with
elementary moves $\{(-\square,\varnothing),(\varnothing,+\square),
(\varnothing,\varnothing),(-\square,+\square)\}$ and partitions and
finally the bijection between the vacillating tableaux with
elementary moves
$\{(-\square,\varnothing),(\varnothing,+\square),(\varnothing,\varnothing),
(+\square,-\square)\}$ and enhanced partitions.


\section{Tangled-diagrams and vacillating tableaux}


\subsection{Tangled-diagrams}
A tangled-diagram is a labeled  graph, $G_n$, over $[n]$ with degree
$\le 2$, represented by drawing its vertices in a
horizontal line and its arcs $(i,j)$ in the upper halfplane having the
following properties:
two arcs $(i_1,j_1)$ and $(i_2,j_2)$ such that $i_1<i_2<j_1<j_2$ are crossing
and if $i_1<i_2<j_2<j_1$ they are nesting.
Two arcs $(i,j_1)$ and $(i,j_2)$ (common lefthand endpoint) and
$j_1<j_2$ can be drawn in two ways: either draw $(i,j_1)$ strictly below
$(i,j_2)$ in which case $(i,j_1)$ and $(i,j_2)$ are nesting (at $i$) or
draw $(i,j_1)$ starting above $i$ and intersecting $(i,j_2)$ once,
in which case $(i,j_1)$ and $(i,j_2)$ are crossing (at $i$):
\begin{figure}[ht]
\epsfig{file=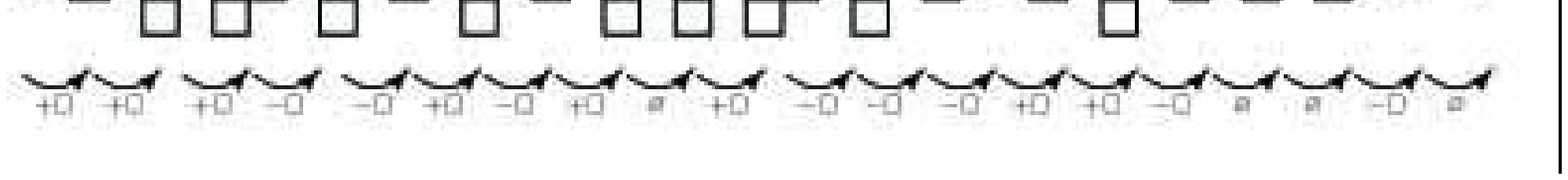,width=0.8\textwidth}
\end{figure}
The cases of two arcs $(i_1,j)$, $(i_2,j)$, where $i_1<i_2$ (common
righthand endpoint)
\begin{center}
\scalebox{0.6}[0.6]{\includegraphics*[60,780][560,830]{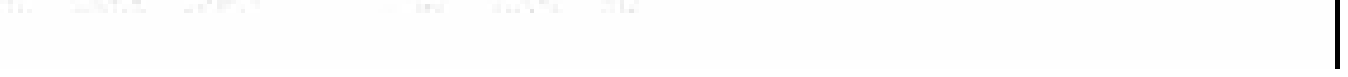}}
\end{center}
and of two arcs $(i,j),(i,j)$, i.e.~where $i$ and $j$ are both:
right- and lefthand endpoints are completely analogous. Suppose
$i<j<h$ and that we are given two arcs $(i,j)$ and $(j,h)$. Then we
can draw them intersecting once or not. In the former case $(i,j)$
and $(j,h)$ are crossing:
\begin{figure}[ht]
\epsfig{file=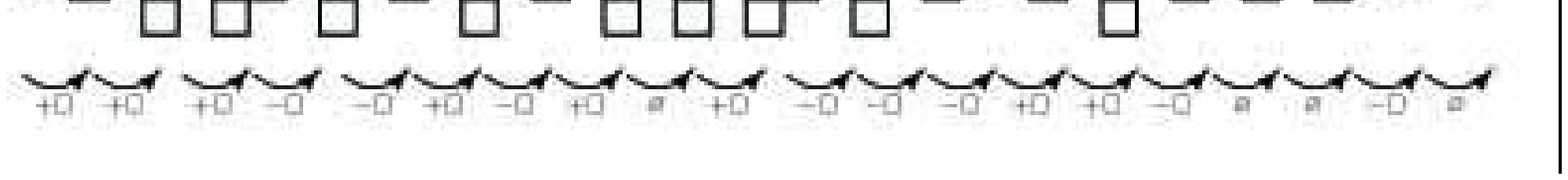,width=0.8\textwidth}
\end{figure}
The cases of two arcs $(i_1,j)$, $(i_2,j)$, where $i_1<i_2$ (common
righthand endpoint)
\begin{center}
\scalebox{0.6}[0.6]{\includegraphics*[60,780][560,830]{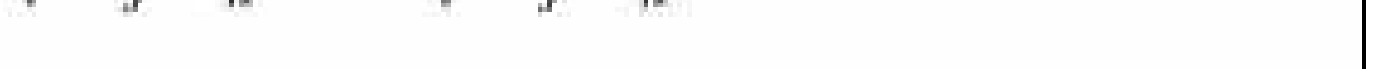}}
\end{center}
The set of all tangled-diagrams is denoted by $\mathcal{G}_n$. A
tangled-diagram is $k$-noncrossing if it contains no $k$-set of
mutually intersecting arcs and $k$-nonnesting if it contains no
$k$-set of mutually nesting arcs. The set of $k$-noncrossing and
$k$-nonnesting tangled-diagrams is denoted by $\mathcal{G}_{n,k}$
and $\mathcal{G}_n^k$, respectively.

\subsection{Inflation}
A key observation allowing for the combinatorial interpretation of
tangled-diagrams is their ``local'' inflation. Intuitively, a tangled-diagram
with $\ell$ vertices of degree $2$ is resolved into a partial matching
over $n+\ell$ vertices. For this purpose we consider the following linear
ordering over
$\{1,1',\dots,n,n'\}$
\begin{equation}
1<1'<2<2'<\dots <(n-1)<(n-1)'<n<n' \ .
\end{equation}
Let $G_n$ be a tangled-diagram with exactly $\ell$ vertices of degree $2$.
Then the inflation of $G_n$, $\iota(G_n)$, is a combinatorial graph over
$\{1,\dots,n+\ell\}$ vertices with degree $\le 1$
obtained as follows:\\
\underline{$i<j_1<j_2$}: if $(i,j_1)$, $(i,j_2)$ are crossing, then
$((i,j_1), (i,j_2))\mapsto ((i,j_1),(i',j_2))$ and
if $(i,j_1)$, $(i,j_2)$ are nesting then
$((i,j_1), (i,j_2))\mapsto ((i,j_2),(i',j_1))$, i.e.:
\begin{figure}[ht]
\epsfig{file=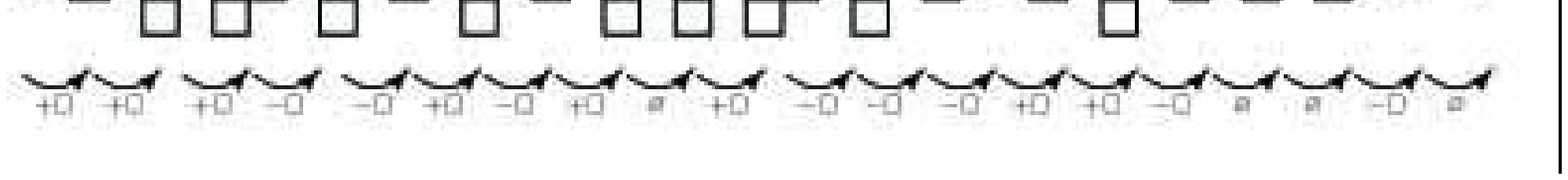,width=0.8\textwidth}
\end{figure}
The cases of two arcs $(i_1,j)$, $(i_2,j)$, where $i_1<i_2$ (common
righthand endpoint)
\begin{center}
\scalebox{0.6}[0.6]{\includegraphics*[60,780][560,830]{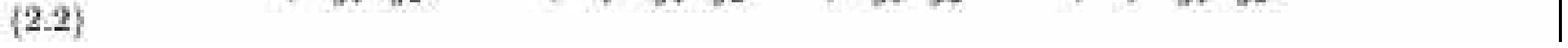}}
\end{center}
\underline{$i_1<i_2<j$}: if $(i_1,j)$, $(i_2,j)$ are crossing then
$((i_1,j), (i_2,j))\mapsto ((i_1,j),(i_2,j'))$ and if $(i_1,j)$,
$(i_2,j)$ are nesting then
$((i_1,j), (i_2,j))\mapsto ((i_1,j'),(i_2,j))$, i.e.:\\
\begin{center}
\scalebox{0.4}[0.5]{\includegraphics*[20,760][560,830]{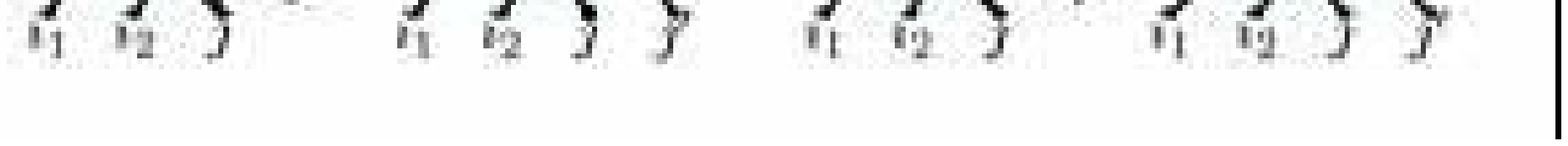}}
\end{center}
\underline{$i<j$}:
if $(i,j),(i,j)$ are crossing, then $((i,j), (i,j))\mapsto ((i,j),
(i',j'))$ and if $(i,j),(i,j)$ are nesting,
then we set $((i,j), (i,j))\mapsto ((i,j'),(i',j))$ and if
$(i,i)$ is a loop we map $(i,i)\mapsto (i,i')$:\\
\begin{center}
\scalebox{0.5}[0.5]{\includegraphics*[20,770][560,820]{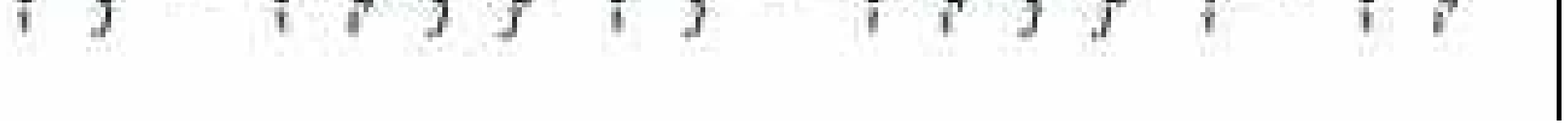}}
\end{center}
\underline{$i<j<h$}: if $(i,j)$, $(j,h)$ are crossing, then $((i,j),
(j,h))\mapsto ((i,j'),(j',h))$ and $((i,j), (j,h))\mapsto ((i,j),(j',h))$,
otherwise, i.e.~we have the following situation
\begin{center}
\scalebox{0.4}[0.5]{\includegraphics*[20,760][560,830]{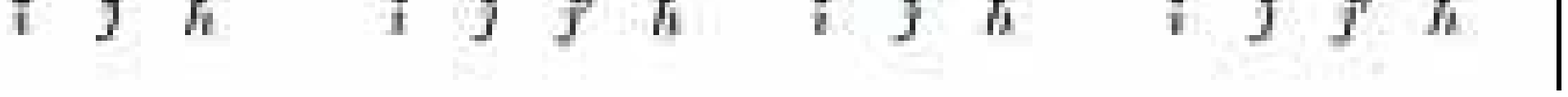}}
\end{center}
Identifying all vertex-pairs $(i,i')$ allows us to recover the original
tangled-diagram and we accordingly have the bijection
\begin{equation}
\iota\colon \mathcal{G}_n \longrightarrow \iota(\mathcal{G}_n) \ .
\end{equation}
$\iota$ preserves by definition the maximal number crossing and nesting arcs,
respectively. Equivalently, a tangled-diagram $G_n$ is $k$-noncrossing if and
only if its inflation $\iota(G_n)$ is $k$-noncrossing or $k$-nonnesting,
respectively. For instance the inflation of the tangled-diagram of
Section~\ref{S:intro} is
\begin{center}
\scalebox{0.6}[0.6]{\includegraphics*[20,770][560,820]{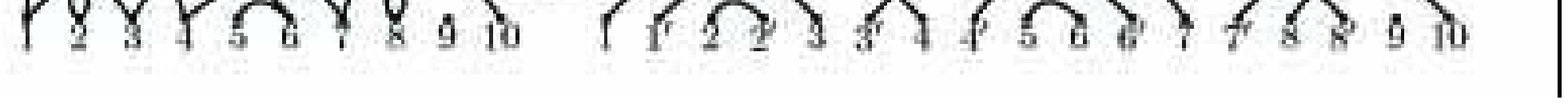}}
\end{center}
\subsection{Vacillating tableaux}
A Young diagram (shape) is a collection of squares arranged in
left-justified rows with weakly decreasing number of boxes in each
row. A standard Young tableau (SYT) is a filling of the squares by
numbers which is strictly decreasing in each row and in each column.
We refer to standard Young tableaux as Young tableaux. Elements can
be inserted into SYT via the RSK-algorithm \cite{ec1}. In the
following we will refer to SYT simply as tableaux. Our first lemma
is instrumental for constructing the bijection between
vacillating tableaux and tangled-diagrams in Section~\ref{S:biject}.

\begin{lemma}{\label{L:extract}}\cite{Chen}
Suppose we are given two shapes $\lambda^{i}\subsetneq\lambda^{i-1}$, which
differ by exactly one square and $T_{i-1}$ is a SYT of shape
$\lambda^{i-1}$. Then there exists a unique $j$ contained in $T_{i-1}$
and a unique tableau $T_{i}$ such that $T_{i-1}$ is obtained from $T_{i}$
by inserting $j$ via the $\text{\rm RSK}$-algorithm.
\end{lemma}


\begin{proof}
Suppose we have two shapes $\lambda^{i}\subsetneq \lambda^{i-1}$, which
differ by exactly one square and $T_{i-1}$ is a tableau of shape
$\lambda^{i-1}$. Let us first assume that $\lambda^{i-1}$ differs from
$\lambda^{i}$ by the rightmost square in its first row, containing $j$.
Then $j$ is the unique element of $T_{i-1}$ which, if $\text{RSK}$-inserted
into $T_i$, produces the tableau $T_{i-1}$.
Suppose next the square which is being removed from $\lambda^{i-1}$
is at the end of row $\ell$. Then we remove the square and RSK-insert its
element $x$ into the $(\ell-1)$-th row in the square containing $y$,
where $y$ is maximal subject to $y< x$ and such that $y$, if inserted
into row $(\ell-1)$, would push down $x$ in its original position. Since
each column is strictly increasing such an $y$ always
exists. We can conclude by induction on $\ell$ that this process results
in exactly one element $j$ being removed from $T_{i-1}$ and a filling of
$\lambda^{i}$, i.e.~a unique tableau $T_i$. By construction,
$\text{\rm RSK}$-insertion of $j$ recovers the tableaux $T_{i-1}$.
\end{proof}

\begin{definition} (Vacillating Tableau)
A vacillating tableaux $V_{\lambda}^{2n}$ of shape $\lambda$ and
length $2n$ is a sequence $(\lambda^{0}, \lambda^{1},\ldots,
\lambda^{2n})$ of shapes such that
{\sf (i)} $\lambda^{0}=\varnothing$ and $\lambda^{2n}=\lambda,$ and
{\sf (ii)} $(\lambda^{2i-1},\lambda^{2i})$ is derived from
           $\lambda^{2i-2}$, for $1\le i\le n$ by either
$(\varnothing,\varnothing)$: do nothing twice;
$(-\square,\varnothing)$: first remove a square then do nothing;
$(\varnothing,+\square)$: first do nothing then adding a square;
$(\pm \square,\pm \square)$: add/remove a square at the odd and even steps,
respectively.
Let $\mathcal{V}_\lambda^{2n}$ denote the set of vacillating
tableaux.
\end{definition}



\section{The bijection}\label{S:biject}


\begin{lemma}\label{L:v-g}
There exists a mapping from the set of vacillating tableaux of shape
$\varnothing$ and length $2n$, $V_\varnothing^{2n}$, into the set of
inflations of tangled-diagrams
\begin{equation}
\phi\colon V_{\varnothing}^{2n}  \longrightarrow \iota(\mathcal{G}_n) \ .
\end{equation}
\end{lemma}

\begin{proof}
In order to define $\phi$ we recursively define a sequence of
triples
\begin{equation}
((P_{0},T_{0},V_0),(P_{1},T_{1},V_{1}),\ldots,(P_{2n},T_{2n},V_{2n}))
\end{equation}
where $P_{i}$ is a set of arcs, $T_{i}$ is a tableau of shape
$\lambda^{i}$, and  $V_i\subset \{1,1',2,2',\dots,n,n'\}$ is a set
of vertices. $P_{0}=\varnothing$, $T_{0}=\varnothing$ and
$V_{0}=\varnothing$. We assume that the lefthand and righthand
endpoints of all $P_{i}$-arcs and the entries of the tableau $T_{i}$
are contained in $\{1,1',\dots,n,n'\}$.
Given $(P_{2j-2},T_{2j-2},V_{2j-2})$ we derive $(P_{2j-1},T_{2j-1},V_{2j-1})$  and $(P_{2j},T_{2j},V_{2j})$ as follows: \\
$\underline{(\varnothing,\varnothing)}$.
If $\lambda^{2j-1}=\lambda^{2j-2}$ and
$\lambda^{2j}=\lambda^{2j-1}$, we have
$(P_{2j-1},T_{2j-1})=(P_{2j-2},T_{2j-2})$ and
$(P_{2j},T_{2j})=(P_{2j-1},T_{2j-1})$ and $V_{2j-1}=V_{2j-2}\cup \{j\}$,
$V_{2j}=V_{2j-1}$.\\
$\underline{(-\square,\varnothing)}$.
If $\lambda^{2j-1}\subsetneq\lambda^{2j-2}$ and
$\lambda^{2j}=\lambda^{2j-1}$, then $T_{2j-1}$ is the unique
tableau of shape $\lambda^{2j-1}$ such that $T_{2j-2}$ is obtained by
RSK-inserting the unique number $i$ via the $\text{RSK}$-algorithm into
$T_{2j-1}$ (Lemma \ref{L:extract}) and
$P_{2j-1}=P_{2j-2}\,\cup\,\{(i,j)\}$ and
$(P_{2j},T_{2j})=(P_{2j-1},T_{2j-1})$ and $V_{2j-1}=V_{2j-2}\cup \{j\}$,
$V_{2j}=V_{2j-1}$. \\
$\underline{(\varnothing, +\square)}$. If $\lambda^{2j-1}=\lambda^{2j-2}$ and
$\lambda^{2j}\supsetneq\lambda^{2j-1}$, then
$(P_{2j-1},T_{2j-1})=(P_{2j-2},T_{2j-2})$ and $P_{2j}=P_{2j-1}$
and $T_{2j}$ is obtained from $T_{2j-1}$ by adding the entry $j$ in the
square $\lambda^{2j}\setminus \lambda^{2j-1}$ and
$V_{2j-1}=V_{2j-2}$, $V_{2j}=V_{2j-1}\cup \{j\}$.\\
$\underline{(-\square,+\square)}$.
If $\lambda^{2j-1}\subsetneq\lambda^{2j-2}$ and
$\lambda^{2j}\supsetneq \lambda^{2j-1}$, then $T_{2j-1}$ is the unique
tableau of shape $\lambda^{2j-1}$ such that $T_{2j-2}$ is obtained
from $T_{2j-1}$ by RSK-inserting the unique number $i$, via
the $\text{RSK}$-algorithm (Lemma \ref{L:extract}). Then
we set
$P_{2j-1}=P_{2j-2}\,\cup\, \{(i,j)\}$ and $P_{2j}=P_{2j-1}$
and $T_{2j}$ is obtained from $T_{2j-1}$ by adding the entry $j'$ in the
square $\lambda^{2j}\setminus \lambda^{2j-1}$ and $V_{2j-1}=
V_{2j-2}\cup\{j\}$, $V_{2j}=V_{2j-1}\cup \{j'\}$. \\
$\underline{(+\square,-\square)}$.
If $\lambda^{2j-2}\subsetneq \lambda^{2j-1}$ and
$\lambda^{2j}\subsetneq\lambda^{2j-1}$ then $T_{2j-1}$ is obtained from
$T_{2j-2}$ by adding the entry $j$ in the square $\lambda^{2j-1}\setminus
\lambda^{2j-2}$ and the tableau $T_{2j}$ is the unique
tableau of shape $\lambda^{2j}$ such that $T_{2j-1}$ is obtained
from $T_{2j}$ by RSK-inserting the unique number $i$, via
the $\text{RSK}$-algorithm (Lemma \ref{L:extract}).
We then set
$P_{2j-1}=P_{2j-2}$ and $P_{2j}=P_{2j-1}\,\cup\, \{(i,j')\}$
and $V_{2j-1}=V_{2j-2}\cup\{j\}$, $V_{2j}=V_{2j-1}\cup \{j'\}$.\\
$\underline{(-\square,-\square)}$.
If $\lambda^{2j-1}\subsetneq \lambda^{2j-2}$ and
$\lambda^{2j}\subsetneq \lambda^{2j-1}$, let $T_{2j-1}$ be the unique
tableau of shape $\lambda^{2j-1}$ such that $T_{2j-2}$ is obtained from
$T_{2j-1}$ by RSK-inserting $i_1$ (Lemma \ref{L:extract}) and $T_{2j}$ be
the unique tableau of shape $\lambda^{2j}$ such that
$T_{2j-1}$ is obtained from $T_{2j}$ by
RSK-inserting $i_2$ (Lemma \ref{L:extract})
$P_{2j-1}=P_{2j-2}\,\cup \, \{(i_1,j)\}$ and
$P_{2j}=P_{2j-1}\,\cup \, \{(i_2,j')\}$ and
$V_{2j-1}=V_{2j-2}\cup\{j\}$, $V_{2j}=V_{2j-1}\cup \{j'\}$.\\
$\underline{(+\square,+\square)}$.
If $\lambda^{2j-1}\supsetneq\lambda^{2j-2}$ and
$\lambda^{2j}\supsetneq \lambda^{2j-1}$, we set
$P_{2j-1}=P_{2j-2}$, and $T_{2j-1}$ is obtained from $T_{2j-2}$ by adding
the entry $j$ in the square $\lambda^{2j-1}\setminus \lambda^{2j-2}$.
Furthermore we set $P_{2j}=P_{2j-1}$ and $T_{2j}$ is obtained from
$T_{2j-1}$ by adding the entry $j'$ in the square $\lambda^{2j}\setminus
\lambda^{2j-1}$ and $V_{2j-1}=V_{2j-2}\cup\{j\}$, $V_{2j}=V_{2j-1}\cup
\{j'\}$.\\
{\it Claim.} $\phi(V_\varnothing^{2n})$ is the inflation of a tangled-diagram.\\
First, if $(i,j)\in P_{2n}$, then $i<j$ and
secondly any vertex $j$ can occur only as either as lefthand or righthand
endpoint of an arc, whence $\phi(V_\varnothing^{2n})$ is a $1$-diagram.
Each step $(+\square,+\square)$ induces a pair of arcs of the form
$(i,j_1)$, $(i',j_2)$ and each step $(-\square,-\square)$ induces a pair
of arcs of the form $(i_1,j)$, $(i_2,j')$.
Each step $(-\square,+\square)$ corresponds to a pair of arcs $(h,j)$,
$(j',s)$ where $h<j<j'<s$ and each step $(+\square,-\square)$ induces
a pair of arcs of the form $(j,s)$, $(h,j')$, where $h<j<j'<s$ or
a $1$-arcs of the form $(i,i')$.
Let $\ell$ be the number of steps not containing $\varnothing$. By
construction each of these adds the $2$-set $\{j,j'\}$, whence
$(V_{2n},P_{2n})$ corresponds to the inflation of a unique tangled-diagram
with $\ell$ vertices of degree $2$ and the claim follows.
\end{proof}
\begin{remark}
The mapping $\phi$: if squares are added the corresponding numbers are
inserted, if squares are deleted Lemma~\ref{L:extract} is used to extract
a unique number, which then forms the lefthand endpoint of the derived arcs.
\begin{center}
\scalebox{0.6}[0.6]{\includegraphics*[20,720][560,820]{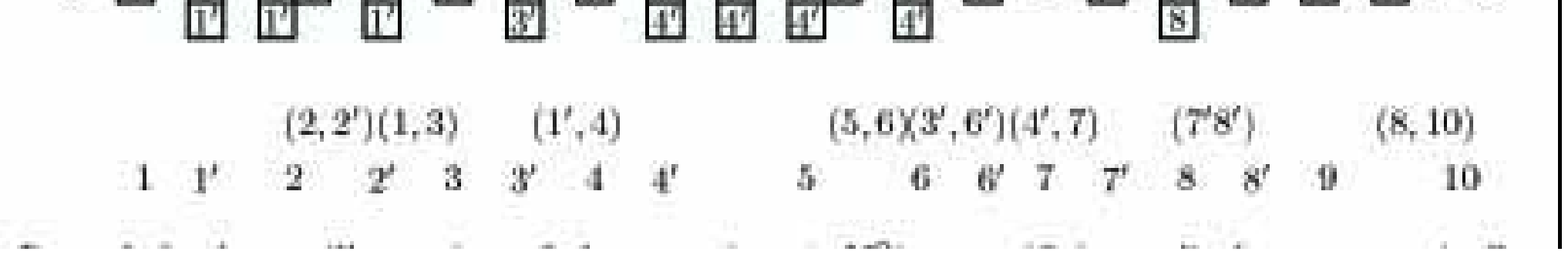}}
\end{center}
\end{remark}
\begin{remark}
As an illustration of the mapping $\phi\colon V_{\varnothing}^{2n}
\longrightarrow \iota(\mathcal{G}_n)$ we display systematically all
arc-configurations of inflated tangled-diagrams induced by the
vacillating tableaux
\begin{center}
\scalebox{0.5}[0.5]{\includegraphics*[20,460][570,820]{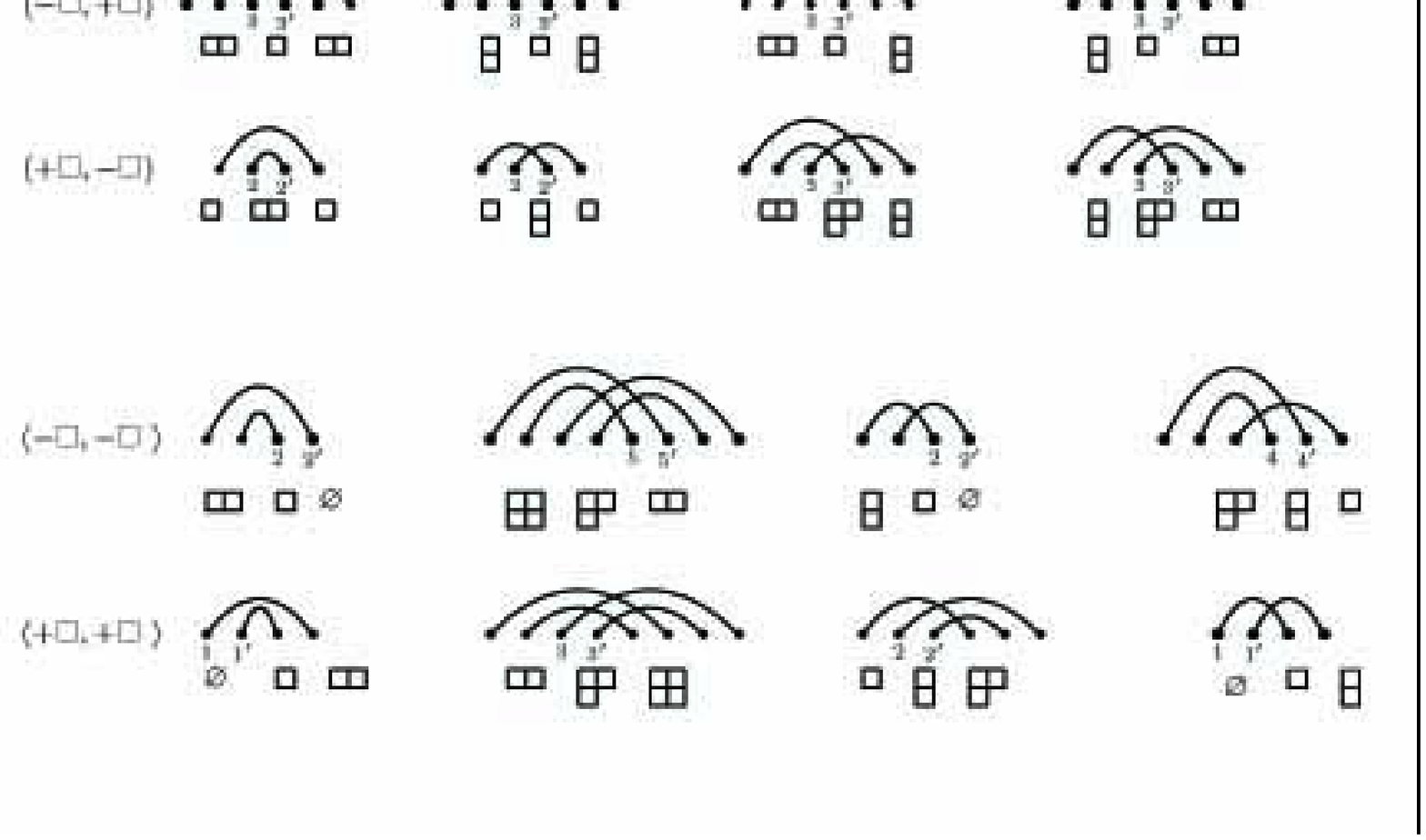}}
\end{center}
\end{remark}

We proceed by explicitly constructing the inverse of $\phi$.

\begin{lemma}\label{L:g-v}
There exists a mapping from the set of inflations of tangled-diagrams
over $n$ vertices, $\iota(\mathcal{G}_n)$, into the set of vacillating
tableaux of shape $\varnothing$ and length $2n$, $\mathcal{V}_{\varnothing}
^{2n}$
\begin{equation}
\psi\colon \iota(\mathcal{G}_n)\longrightarrow
\mathcal{V}_{\varnothing}^{2n} \ .
\end{equation}
\end{lemma}
\begin{proof}
We define $\psi$ as follows. Let $\iota(G_n)$ be the inflation of the
tangled-diagram $G_n$. We set
\begin{equation}\label{E:iota}
\iota_i=
\begin{cases}
(i,i') & \text{\rm iff } \ \text{\rm $i$ has degree $2$ in $G_n$}, \\
i      & \text{\rm otherwise.}
\end{cases}
\end{equation}
Let $T_{2n}=\varnothing$ be the empty tableau. We will construct inductively
a sequence of tableaux $T_h$ of shape $\lambda_{\iota(G_n)}^{h}$, where $h\in
\{0,1,\ldots 2n\}$ by considering $\iota_i$ for $i=n,n-1,n-2,\dots, 1$.
For each $\iota_j$ we inductively define the pair of tableaux $(T_{2j},
T_{2j-1})$:\\
{\sf (I)} $\iota_j=j$ is an isolated vertex in $\iota(G_n)$, then we
set $T_{2j-1}=T_{2j}$ and $T_{2j-2}=T_{2j-1}$.
Accordingly, $\lambda_{\iota(G_n)}^{2j-1}=
\lambda_{\iota(G_n)}^{2j}$ and $\lambda_{\iota(G_n)}^{2j-2}=
\lambda_{\iota(G_n)}^{2j-1}$ (left to right: $(\varnothing,\varnothing)$).\\
{\sf (II)} $\iota_j=j$ is the righthand endpoint of exactly one arc
$(i,j)$ but not a lefthand endpoint, then we set $T_{2j-1}=T_{2j}$ and obtain
$T_{2j-2}$ by adding $i$ via the $\text{\rm RSK}$-algorithm to
$T_{2j-1}$. Consequently we have $\lambda_{\iota(G_n)}^{2j-1}=
\lambda_{\iota(G_n)}^{2j}$ and $\lambda_{\iota(G_n)}^{2j-2} \supsetneq
\lambda_{\iota(G_n)}^{2j-1}$. (left to right: $(-\square,\varnothing)$).\\
{\sf (III)} $j$ is the lefthand endpoint of exactly one arc
$(j,k)$ but not a righthand endpoint, then first set $T_{2j-1}$
to be the tableau obtained by removing the square with entry $j$ from $T_{2j}$
and let $T_{2j-2}=T_{2j-1}$. Therefore $\lambda_{\iota(G_n)}^{2j-1}
\subsetneq\lambda_{\iota(G_n)}^{2j}$ and $\lambda_{\iota(G_n)}^{2j-2}=
\lambda_{\iota(G_n)}^{2j-1}.$ (left to right: $(\varnothing,+\square)$).\\
{\sf (IV)} $j$ is a lefthand and righthand endpoint, then
we have the two $\iota(G_n)$-arcs $(i,j)$ and $(j',h)$,
where $i<j<j'<h$. $T_{2j-1}$ is obtained by removing the square with entry
$j'$ in $T_{2j}$ first and $T_{2j-2}$ via inserting $i$ in $T_{2j-1}$ via the
$\text{\rm RSK}$-algorithm. Accordingly we derive the shapes
$\lambda_{\iota(G_n)}^{2j-1}\subsetneq \lambda_{\iota(G_n)}^{2j}$ and
$\lambda_{\iota(G_n)}^{2j-2}\supsetneq\lambda_{\iota(G_n)}^{2j-1}$.
(left to right: $(-\square,+\square)$). \\
{\sf (V)} $j$ is a righthand endpoint of degree $2$, then we have
the two $\iota(G_n)$-arcs $(i,j)$ and $(h,j')$.
$T_{2j-1}$ is obtained by inserting $h$ via the
$\text{\rm RSK}$-algorithm into $T_{2j}$ and $T_{2j-2}$ is obtained by
RSK-inserting $i$ into $T_{2j-1}$ via the $\text{\rm RSK}$-algorithm.
We derive
$\lambda_{\iota(G_n)}^{2j-1}\supsetneq\lambda_{\iota(G_n)}^{2j}$ and
$\lambda_{\iota(G_n)}^{2j-2}\supsetneq\lambda_{\iota(G_n)}^{2j-1}$ (left to
right:$(-\square,-\square)$). \\
{\sf (VI)} $j$ is a lefthand endpoint of degree $2$, then we have
the two $\iota(G_n)$-arcs $(j,r)$ and $(j',h)$. $T_{2j-1}$ is obtained
by removing the square with entry $j'$ from the tableau $T_{2j}$ and
$T_{2j-2}$ is obtained by removing the square with entry $j$ from the
$T_{2j-1}$.
Then we have $\lambda_{\iota(G_n)}^{2j-1}\subsetneq\lambda_{\iota(G_n)}^{2j}$
and $\lambda_{\iota(G_n)}^{2j-2}\subsetneq\lambda_{\iota(G_n)}^{2j-1}$
(left to right: $(+\square,+\square)$).\\
{\sf (VII)} $j$ is a lefthand and righthand endpoint of crossing arcs
or a loop,
then we have the two $\iota(G_n)$-arcs $(j,s)$ and $(h,j')$, $h<j<j'<s$ or
an arc of the form $(j,j')$.
$T_{2j-1}$ is obtained by RSK-inserting $h$ ($j$) into the tableau
$T_{2j}$ and $T_{2j-2}$ is obtained by removing the square with entry
$j$ ($j$) from the $T_{2j-1}$ (left to right: $(+\square,-\square)$).\\
Therefore $\psi$ maps the inflation of a tangled-diagram into a vacillating
tableau and the lemma follows.
\end{proof}
\begin{remark}
From inflations of tangled-diagrams to vacillating tableaux: starting from
right to left the vacillating tableaux is obtained via the RSK-algorithm
as follows: if $j$ is a righthand endpoint it gives rise to RSK-insertion
of its (unique) lefthand endpoint and if $j$ is a lefthand endpoint the
square containing $j$ is removed.
\begin{center}
\scalebox{0.5}[0.5]{\includegraphics*[20,660][560,820]{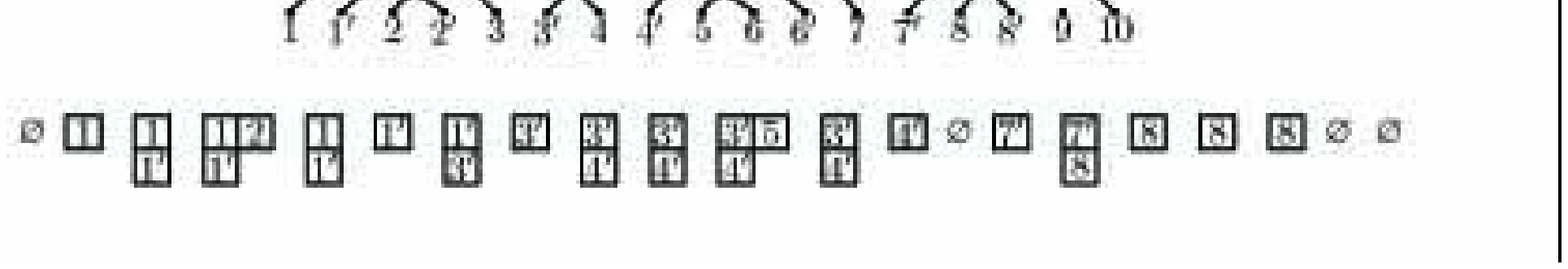}}
\end{center}
\end{remark}
\begin{theorem}\label{T:bij}
There exists a bijection between the set of vacillating tableaux of shape
$\varnothing$ and length $2n$, $\mathcal{V}_\varnothing^{2n}$ and the
set of tangled-diagrams over $n$ vertices, $\mathcal{G}_n$
\begin{equation}
\beta\colon \mathcal{V}_{\varnothing}^{2n}  \longrightarrow
\mathcal{G}_n \ .
\end{equation}
\end{theorem}
\begin{proof}
According to Lemma~\ref{L:v-g} and Lemma~\ref{L:g-v} we have the following
mappings
$\phi\colon \mathcal{V}_{\varnothing}^{2n}  \longrightarrow
\iota(\mathcal{G}_n)$ and
$\psi\colon \iota(\mathcal{G}_n)  \longrightarrow
\mathcal{V}_{\varnothing}^{2n}$.
We next show that $\phi$ and $\psi$ are indeed inverses with respect to each
other. By definition $\phi$ extracts arcs such that their respective
lefthand-endpoints if RSK-inserted (Lemma~\ref{L:extract}) recover the
tableaux of the preceding step. We observe that by definition, $\psi$
reverses this extraction: it explicitly RSK-inserts the lefthand-endpoints
of arcs. Therefore we have the following situation
\begin{equation}
\phi\circ\psi(\iota(G_n))=\phi((\lambda_{\iota(G_n)})_{0}^{2n})=\iota(G_n)
\qquad \text{\rm and}
\qquad \psi\circ \phi(\mathcal{V}_\varnothing^{2n})=
\mathcal{V}_\varnothing^{2n} \ ,
\end{equation}
from which we conclude that $\phi$ and $\psi$ are bijective. Since $G_n$
is in one to one correspondence with $\iota(G_n)$ the proof of the theorem is
complete.
\end{proof}
By construction the bijection $\iota\colon \mathcal{G}_n
\longrightarrow \iota(\mathcal{G}_n)$ preserves the maximal number
crossing and nesting arcs, respectively. Equivalently, a tangled-diagram
$G_n$ is $k$-noncrossing if and only if its inflation $\iota(G_n)$
is $k$-noncrossing or $k$-nonnesting \cite{Chen}. Indeed, this
follows immediately from the definition of the inflation.

\begin{theorem}\label{T:cross}
A tangled-diagram $G_n$ is $k$-noncrossing if and only if all shapes $\lambda^i$
in its corresponding vacillating tableau have less than $k$ rows,
i.e.~$\phi\colon \mathcal{V}_\varnothing^{2n}\longrightarrow \mathcal{G}_n$
maps vacillating tableaux having less than $k$ rows into $k$-noncrossing
tangled-diagrams. Furthermore there is a bijection between the set of
$k$-noncrossing and $k$-nonnesting tangled-diagrams.
\end{theorem}
Theorem~\ref{T:cross} is the generalization of the corresponding result
in \cite{Chen} to tangled-diagrams. Since the inflation map allows to interpret
a tangled-diagram with $\ell$ vertices of degree $2$ over $n$ vertices as a
partial matching over $n+\ell$ vertices its proof is analogous.

We next observe that restricting the steps for vacillating tableaux
produces the bijections of Chen {\it et.al.} \cite{Chen}.
Let $\mathcal{M}_k(n)$, $\mathcal{P}_k(n)$ and $\mathcal{B}_k(n)$
denote the set
of $k$-noncrossing matchings, partitions and braids, respectively.
If a tableaux-sequence $V_\varnothing^{2n}$ is obtained via certain
steps $s\in S$ we write $V_\varnothing^{2n}\models S$.
\begin{corollary}\label{C:restrict}
Let $\beta_i$ denote the restriction of the bijection
$\beta\colon \mathcal{V}_\varnothing^{2n}\longrightarrow \mathcal{G}_n$
in Theorem~\ref{T:bij}. Then $\beta$ induces the bijections
\begin{equation}
\beta_1\colon \{V_\varnothing^{2n}\mid V_\varnothing^{2n} \models
(-\square,\varnothing),(\varnothing,+\square)\, \,\text{\rm  and has
$\le k$ rows} \}  \rightarrow  \mathcal{M}_k(n) \ .
\end{equation}
\begin{equation}
\beta_2\colon \{V_\varnothing^{2n}\mid V_\varnothing^{2n} \models
(-\square,\varnothing),
          (\varnothing,+\square), (\varnothing,\varnothing),
          (-\square,+\square)\,\,\text{\rm  and has $\le k$
rows} \} \rightarrow  \mathcal{P}_k(n) \ .
\end{equation}
\begin{equation}
\beta_3\colon \{V_\varnothing^{2n}\mid V_\varnothing^{2n} \models
(-\square,\varnothing),(\varnothing,+\square),(\varnothing,\varnothing),
(+\square,-\square) \, \,\text{\rm  and has $\le k$ rows} \}
\rightarrow  \mathcal{B}_k(n) \ .
\end{equation}
\end{corollary}
\begin{remark}
For partitions we can illustrate the correspondences between the elementary
steps and associated tangled-diagram arc-configurations as follows:
\begin{center}
\scalebox{0.5}[0.5]{\includegraphics*[30,760][560,820]{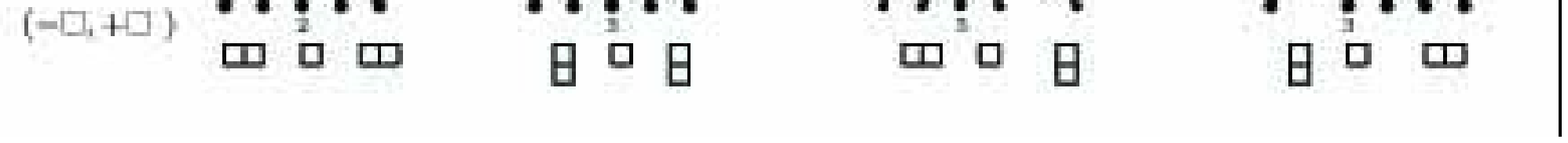}}
\end{center}
\end{remark}
\begin{remark}
For braids we derive the following correspondences. They illustrate
one key difference between partitions and braids: for fixed crossing
number braids are more restricted due to the fact that they already
have {\it a priori} ``local'' crossings at their non-loop-vertices
of degree $2$.
\begin{center}
\scalebox{0.5}[0.5]{\includegraphics*[40,760][560,820]{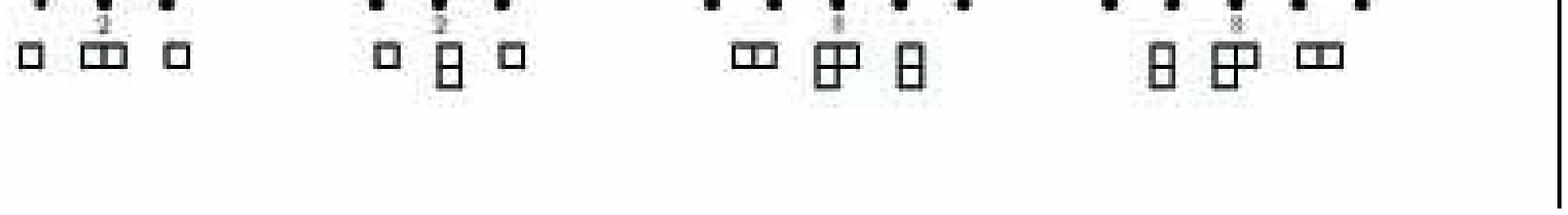}}
\end{center}
\end{remark}

Let $D_{2,k}(n)$ and $\tilde{D}_{2,k}(n)$ be the numbers of
$k$-noncrossing tangled-diagrams and tangled-diagrams without isolated
points over $[n]$, respectively. Furthermore let $f_k(2n-\ell)$ be
the number of $k$-noncrossing matchings over $2n-\ell$ vertices. We
show that the enumeration of tangled-diagrams can be reduced to the
enumeration of matchings via the inflation map. W.l.o.g.~we can
restrict our analysis to the case of tangled-diagrams without isolated
points since the number of tangled-diagrams over $[n]$ is given by
$D_{2,k}(n)=\sum_{i=0}^n\binom{n}{i} \tilde{D}_{2,k}(n-i)$.

\begin{theorem}{\label{T:2diag}}
The number of $k$-noncrossing tangled-diagrams over $[n]$ is
given by
\begin{equation}
\tilde{D}_{2,k}(n)=\sum_{\ell=0}^n\binom{n}{\ell}f_{k}(2n-\ell) \ .
\end{equation}
and in particular for $k=3$ we have
\begin{equation}
\tilde{D}_{2,3}(n)=\sum_{\ell=0}^n\binom{n}{\ell} \left(C_{\frac{2n-\ell}{2}}\,
C_{\frac{2n-\ell}{2}+2}-
C_{\frac{2n-\ell}{2}+1}^2\right)
\ .\\
\end{equation}
\end{theorem}
\begin{proof}
Let $\tilde{\mathcal{D}}_{2,k}(n,V)$ denote the set of tangled-diagrams in
which $V=\{i_1,\dots,i_h\}$
is the set of vertices of degree $1$ (where $h\equiv 0\mod 2$ by definition
of $\tilde{\mathcal{D}}_{2,k}(n,V)$)
and let
$\mathcal{M}_k(\{1,1',\dots,n,n'\}\setminus V')$, where
$V'=\{i_1',\dots,i_h'\}$ denote the set of matchings
over $\{1,1',\dots,n,n'\}\setminus V'$.
By construction, (eq.~(\ref{E:inf1}), eq.~(\ref{E:inf2}) and
eq.~(\ref{E:inf3})) the inflation is a well defined mapping
\begin{equation}
\iota \colon \tilde{\mathcal{D}}_{2,k}(n,V) \longrightarrow
\mathcal{M}_k(\{1,1',\dots,n,n'\}\setminus V')
\end{equation}
with inverse $\kappa$ defined by identifying all pairs $(x,x')$, where
$x,x'\in\{1,1',\dots,n,n'\}\setminus V'$. Obviously, we have
$\vert \mathcal{M}_k(\{1,1',
\dots,n,n'\}\setminus V')\vert=f_k(2n-\ell)$ and we obtain
\begin{equation}
\tilde{D}_{2,k}(n)=\sum_{V\subset [n]}\tilde{D}_{2,k}(n,V)=
\sum_{\ell=0}^n\binom{n}{\ell} f_k(2n-\ell)\ .
\end{equation}
Suppose $n\equiv 0\mod 2$ and let $C_m$ denote the $m$-th Catalan
number, then we have \cite{GB89}
\begin{equation}
f_3(n)=C_{\frac{n}{2}}\,
C_{\frac{n}{2}+2}-
C_{\frac{n}{2}+1}^2
\ .\\
\end{equation}
and the theorem follows.
\end{proof}

\begin{remark}
The first $10$ numbers of 3-noncrossing tangled-diagrams are given by
\\
$$
\begin{tabular}{|c|c|c|c|c|c|c|c|c|c|c|}
  \hline
  $n$ & 1 & 2 & 3 & 4 & 5 & 6 & 7 & 8 & 9 & 10 \\
  \hline
  $D_{2,3}(n)$ & 2 & 7 & 39 & 292 & 2635 & 27019 & 304162 & 3677313 &
47036624 &
  629772754
\\
  \hline
\end{tabular}
$$
\end{remark}

The enumeration of $3$-noncrossing partitions and $3$-noncrossing enhanced
partitions, which are in bijection with braids without isolated points has
been derived in \cite{MIRXIN}.

{\bf Acknowledgments.}
We are grateful to Emma Y.~Jin for helpful discussions. This work was
supported by the 973 Project, the PCSIRT Project of the Ministry of
Education, the Ministry of Science and Technology, and the National
Science Foundation of China.

\bibliographystyle{amsplain}


\end{document}